\magnification=1200
\hsize=16 true cm    \vsize=22.7 true cm
\hoffset=4 true mm   \voffset=2 true mm
\font\tengoth = eufm10
\font\sevengoth = eufm7
\font\fivegoth = eufm5
\newfam\gothfam
\textfont\gothfam=\tengoth
\scriptfont\gothfam=\sevengoth
\scriptscriptfont\gothfam=\fivegoth
\def\goth{\tengoth\fam\gothfam}
\font\tenmath =msbm10
\font\sevenmath=msbm7

\newfam\mathfam
\textfont\mathfam=\tenmath
\scriptfont\mathfam=\sevenmath
\scriptscriptfont\mathfam=\sevenmath
\def\math{\tenmath\fam\mathfam}
\def\build#1_#2{\mathrel{\mathop{\kern 0pt#1}\limits_{#2}}}
\def\K{{\math K}}
\def\Ker{\hbox
{Ker}}
\def\dim{\hbox{dim}}\def\gwo{\Gamma_{\tilde w_0}}
\def\Gr{\hbox{Gr}}%

\def\dim{\hbox{dim}}
\def\K{{\math K}}

\def\build#1_#2{\mathrel{\mathop{\kern 0pt#1}\limits_{#2}}}
\def\B{{\cal B}}\def\A{{A}}\def\cgu{\C[G/U]}
\def\F{{\cal F}}
\def\qed{\hfill$\diamondsuit$}
\def\co{{\cal C}}

\def\ua{U_{\A}}
\def\uqnm{U(\nm)}

\def\R{{\math R}}
\def\uq{U_q(\g)}\def\uqn{U_q(\n)}\def\uqnm{U_q(\nm)}\def\uo{U_q^0}\def\uqb{U_q(
\b)}\def\uqbm{U_q(\b^-)}\def\uanm{U_A(\n^-)}\def\lal{V_A(\lambda)}
\def\lq{V_q}\def\lql{\lq(\lambda)}
\def\cqg{\C_q[G]}\def\cqgu{\C_q[G/U]}
\def\aqgu{A_q[G/U]}


\def\N{{\math N}}\def\Z{{\math Z}}%
\def\Q{{\math Q}}%
\def\n{{\goth n}}%
\def\C{{\math C}}%
\def\g{{\goth g}}%
\def\h{{\goth h}}%
\def\b{{\goth b}}\def\nm{{\goth n}^-}\def\bm{{\goth b^-}}
%

\def\gw{\Gamma_{\tilde w_0}^w}\def\gwp{\Gamma_{\tilde w_0'}^w}

\font\petittitre=cmbx8 scaled 1000 

\vskip 2cm\centerline{\petittitre{ TORIC DEGENERATIONS OF SCHUBERT
VARIETIES}}
\centerline{\petittitre{ }}
\vskip 1cm\centerline{Philippe CALDERO}
\vskip 4cm\noindent
\par\noindent{\petittitre ABSTRACT.} {\it Let $G$ be a simply connected
semi-simple complex algebraic group. We prove that every Schubert variety of $G$ has
a flat degeneration into a toric variety. This provides a
generalization of results of [7], [6], [5]. Our basic tool is
Lusztig's canonical basis and the string parametrization of this basis.}
  \vskip 2cm\noindent
This research has been partially supported by the EC TMR network
 "Algebraic Lie Representations" , contract no. ERB FMTX-CT97-0100.
\vskip 1cm\noindent
\noindent
{\bf 0. Introduction.}\bigskip\noindent
{\bf 0.1.} Let $G$ be a simply connected semi-simple complex algebraic group. 
Fix a
maximal torus $T$ and a Borel subgroup $B$ such that $T\subset B\subset G$. Let
$W$ the Weyl group of $G$ relative to $T$. For any $w$ in $W$, let
$X_w=\overline {BwB/B}$ denote the Schubert variety corresponding to $w$.
This article is concerned with the following problem.\medskip\noindent
{\bf Degeneration Problem.} {\it Is there a flat family over Spec$\C[t]$, such that the
general fiber is $X_w$ and the special fiber is a toric
variety?}\medskip\noindent
The existence of such a degeneration was obtained by N. Gonciulea and V.
Lakshmibai for $G=SL_n$, [7]. Their proof is based on the theory
of standard monomials. In the case $G=SL_n$, the corner stone of their
proof is the following : fundamental weights are minuscule weights, hence,
a basis of every fundamental representation is endowed with a structure of
distributive lattice. \par
A toric degeneration for Schubert varieties is given in [5], [6], for
$G$ of rank 2. The proofs rely on the theory of standard monomials as well.
A natural question would be : is there a (flat) toric  degeneration of the
flag variety $G/B$ which restricts to a toric degeneration of 
the Schubert varieties $X_w$ for any $w$ in the Weyl group? In
[4], R. Chirivi gives a degeneration of the flag variety which  
restricts into semi-toric
degenerations of the Schubert varieties, i.e. finite unions of irreducible
toric varieties. An explanation of this fact was given to us by O. Mathieu
: intersections of irreducible toric varieties are irreducible toric
varieties, but intersection of Schubert varieties can be a union of several 
Schubert
varieties. Hence, the answer to the previous question is negative. In [4],
the Degeneration Problem is solved with toric replaced by semi-toric. \smallskip\noindent
{\bf 0.2.}\def\ll{{\cal L}_\lambda}\def\lml{{\cal L}_{m\lambda}}
Our approach of the problem is based on the canonical/global base of
Lusztig/Ka\-shi\-wa\-ra and the so-called string parametrization of this base
studied by P. Littelmann in [10] and made precise by A. Berenstein and A.
Zelevinsky in [1].\par
Fix $w$ in $W$. 
Let $P^+$ be the semigroup of dominant weights. For all $\lambda$ in 
$P^+$, let $\ll$ be the line bundle on $G/B$ corresponding to $\lambda$. Then,
the direct sum of global sections $R_w:=\bigoplus_{\lambda\in
P^+}H^0(X_w,\ll)$
carries a natural structure of $P^+$-graded $\C$-algebra. Moreover, there
exists a natural action of $T$ on $R_w$. Our principal result can be
stated as follows :\medskip\noindent
{\bf Theorem.} {\it Fix $w$ in $W$. There exists a filtration
$(\F_m^w)_{m\in\N}$ of $R_w$ such that\smallskip\noindent
\item{(i)} for all $m$ in $\N$, $\F_m^w$ is compatible with the   
$P^+$-grading of $R_w$,\par\noindent
\item{(ii)} for all $m$ in $\N$, $\F_m^w$ is compatible with the action of
$T$, \par\noindent 
\item{(iii)} the associated graded algebra is the $\C$-algebra of the
semigroup of integral points in a rational convex polyhedral cone.}\medskip\noindent
This cone depends on the choice of a reduced decomposition $\tilde w_0$ of
the longest element $w_0$ of the Weyl group. Explicit equations for the
faces of this cone can be obtained from [10] for so-called nice decompositions
$\tilde w_0$. More generally, those equations were obtained in [1]
from $\tilde w_0$-trails in fundamental Weyl modules of the Langlands 
dual of $G$. \par
This theorem gives a positive answer to the Degeneration Problem. Indeed, 
let $\lambda$ be a regular dominant weight, then the line bundle $\ll$ is
ample and $X_w$ is the projective spectrum of $\bigoplus_{m\in
\N}H^0(X_w,\lml)$. Moreover, the spectrum of a
noetherian graded algebra associated to a filtration of a noetherian
algebra $R$ is a flat degeneration of Spec$(R)$. This is proved as follows by 
a standard argument : let $t$ be an indeteminate and consider the filtration
$(R_n)_{n\in\N}$ of $R$. Then, the $\C[t]$-algebra $R^t=\oplus_nR_nt^n$ is flat
over $\C[t]$ and it verifies $R^t/(t-t_0)R^t\simeq R$ for $t_0\not=0$ and
$R^t/tR^t\simeq \Gr R$.\par
Let $U$ be the maximal unipotent subgroup of $B$. Hence, $R_{w_0}$ is the
algebra $\cgu$ of regular functions on $G/U$. For any $w$, the
algebras $R_w$ are quotients of this algebra. If $\tilde w_0$
is adapted to $w$ in the sense of Definition 2.4, then the filtration
$(\F_m^w)_{m\in\N}$ is the quotient filtration of $(\F_m^{w_0})_{m\in\N}$.
In general, the quotient filtration  $(\overline\F_m^{w_0})_{m\in\N}$ of
$R_w$  provides a graded associated algebra whose spectrum is a semi-toric
variety.\smallskip\noindent
The proof of the theorem is based on two facts.\smallskip\noindent  
Let  $U^-$ be the  maximal unipotent subgroup of $G$ which is opposite to $U$.
Set $B^-=TU^-$. Then, the algebra $\cgu$ embeds in $C[B^-]$. Moreover,
 we can embed the (specialized) dual of the canonical base in
the algebra $\C[U^-]$. We prove that this dual has good
multiplicative properties inherited from the quantum case, see Theorem 2.3. This part of the article is
inspired by [13]. But here, we don't use the positivity arguments or the elaborate Hall
algebra model, only true for the simply laced case.
\smallskip\noindent
In a second step we show how to restrict from flag
variety to Schubert varieties : this part relies on the
 compatibility of the canonical base with
the Demazure modules, [9], [10].
\bigskip\noindent
{\bf 1. Notations and recollection on global basis.}\bigskip\noindent
{\bf 1.1.} Denote by $G$ a semisimple simply connected complex Lie group. Fix
a torus $T$ of
$G$ and let $B$ be a
Borel subgroup such that $T\subset B\subset G$. Denote by
 $U$ the unipotent radical of
$B$.
Let $B^-$ be the opposite Borel subgroup and $U^-$
be its unipotent radical. 
Let $\g$, resp. $\h$, $\n$, $\b$, $\nm$, $\bm$, be 
 Lie $\C$-algebra of $G$, resp. $T$, $U$, $B$, $U^-$, $B^-$. Let $n$
be the rank of $\g$. 
We have the triangular decomposition $\g=\nm+\h+\n$.
Let $\{\alpha_i\}_i$ be a basis of the root system
$\Delta$ corresponding to this decomposition. 
Let  $P$ be the weight lattice generated by the fundamental weights
$\varpi_i$, $1\leq i\leq n$, 
and let $P^+:=\sum_i\,\N\varpi_i$ be the semigroup of integral dominant weights.
Let $W$
 be the  Weyl group, generated by the
 reflections $s_{\alpha_i}$ corresponding to the simple roots $\alpha_i$,
 and let $w_0$ be the longest element of $W$.
We denote by $(\,,\,)$ the $W$-invariant form on $P$.
\medskip\noindent
{\bf 1.2.} Let $d$ be an integer such that $(P,P)\subset (2/d)\Z$.
Let $q$ be a indeterminate and set $\K=\C(q^{1/d})$.
Let $\uq$ be the simply connected quantized
enveloping  algebra on $\K$, as defined in [3]. Set
$d_i=(\alpha_i,\alpha_i)/2$
and $q_i=q^{d_i}$ for all $i$. 
Let $\uqn$, resp. $\uqnm$, 
be the subalgebra generated by the canonical generators  
$E_{\alpha_i}$, resp. $F_{\alpha_i}$,
 of positive,  resp. negative, weights and the quantum Serre relations.
For all $\lambda$
in $P$, let $K_\lambda$ the corresponding  element in the 
algebra $\uo=\K[P]$ of the torus of $\uq$.
 We have the triangular decomposition   
$\uq=\uqnm\otimes\uo\otimes\uqn$. We set 
$$\uqb=\uqn\otimes\uo, \hskip 15mm \uqbm=\uqnm\otimes\uo.$$\par
$\uq$ is endowed with a  structure of Hopf algebra and the 
comultiplication $\Delta$, 
the antipode $S$ and the augmentation $\varepsilon$ are given by
$$\Delta E_i=E_i\otimes 1+K_{\alpha_i}\otimes E_i,\,\Delta F_i
=F_i\otimes K_{\alpha_i}^{-1}+ 1\otimes F_i, \Delta K_\lambda
=K_\lambda\otimes K_\lambda$$
$$S(E_i)=-K_{\alpha_i}^{-2}E_i,\,S(F_i)=-F_iK_{\alpha_i}^{2}, S(K_\lambda)
=K_{-\lambda}$$
$$\varepsilon(E_i)=\varepsilon(F_i)=0,\,\varepsilon(K_\lambda)
=1.$$\medskip\noindent
If $M$ is a $\uo$-module and $\gamma\in P$, we set $M_\gamma:=\{m\in M,\;
K_\lambda.m=q^{(\lambda,\gamma)}\}$.\par
For 
$n$  a non negative integer and $\alpha$ a positive root, we set  :
$[n]_i={1-q_i^n\over 1-q_i}$, $[n]_i!=[n]_i
[n-1]_i\ldots[1]_i$.\medskip\noindent
{\bf 1.3.} The dual $\uq^*$ is endowed with a structure of left, 
resp. right, $\uq$-module by 
$u.c(a)=c(au)$, resp. $c.u(a)=c(ua)$, $u,\,a\in\uq$, $c\in\uq^*$. 
If $M$ is a finite dimensional left $\uq$-module, 
we endow the dual $M^*$ with a structure of
left $\uq$-module by $u\xi(v)=\xi(S(u)v)$, $u\in\uq$, $\xi\in M^*$, $v\in
M$.\par
For all $\lambda$ in $P^+$, let $\lql$ be the simple $\uq$-module 
with highest weight $\lambda$. We can embed $\lql^*\otimes\lql$ in $\uq^*$ by
setting $\xi\otimes v(u)=\xi(u.v)$, $u\in\uq$, $\xi\in\lql^*$, $v\in\lql$.   
Let $v_\lambda$ be a highest weight vector of $\lql$.
For all integral dominant weight $\lambda$, let $C(\lambda)$, resp.
$C^+(\lambda)$, 
be the subspace of $\uq^*$ generated by the $\xi\otimes v$, resp.
$\xi\otimes v_\lambda$, $\xi\in\lql^*$, $v\in\lql$. We set 
$\cqg=\bigoplus_{\lambda\in P^+}\,C(\lambda)$, 
$\cqgu=\bigoplus_{\lambda\in P^+}\,C^+(\lambda)$. Then,
$\cqg$ and $\cqgu$  are subalgebras of the Hopf dual of $\uq$. 
$\cqg$, resp. $\cqgu$, is the algebra of quantum regular functions on $G$,
resp. on the quotient $G/U$.\medskip\noindent
{\bf 1.4.} There exists a unique bilinear form $(\,,\, )$ on
$\uqb\times\uqbm$, see [14], [15], [3],
such that :
$$(u^+,\,u_1^-u_2^-)=(\Delta(u^+),\,u_1^-\otimes u_2^-)\,,\hskip 1cm
u^+\in\uqb\,;\,u_1^-,u_2^-\in \uqbm\leqno (1.4.1)$$
$$(u_1^+u_2^+,\,u^-)=(u_2^+\otimes u_1^+,\,\Delta(u^-))\,,\hskip 1cm
u^-\in\uqbm\,;\,u_1^+,u_2^+\in\uqb\leqno (1.4.2)$$
$$(K_\lambda,\,K_\mu)=q^{-(\lambda,\mu)}\,,\hskip 1cm
\lambda,\mu\in P\leqno (1.4.3)$$
$$(K_\lambda,\,F_i)=0\,,\hskip 1cm\lambda\in P\,,\,1\leq i\leq
n\leqno (1.4.4)$$
$$(E_i,\,K_\lambda)=0\,,\hskip 1cm \lambda\in P\,,\,1\leq i\leq
n\leqno (1.4.5)$$
$$(E_i,\,F_i)=\delta_{ij}(1-q_i^2)^{-1}\,,\hskip 1cm 1\leq i,j \leq
n.\leqno (1.4.6)$$
For all $\beta$ in $Q^+$, let $\uqn_\beta$, resp. $\uqnm_{-\beta}$, be the 
subspace of $\uqn$, resp. $\uqnm$, with weight $\beta$, resp. $-\beta$.
The form $(\,,\, )$ is non degenerate on $\uqn_\beta
\times \uqnm_{-\beta}$, $\beta\in Q^+$. We have, by (1.4.1-1.4.5) :
$$(XK_\lambda,\,YK_\mu)=q^{-(\lambda,\mu)}(X,\,Y)\,,
\hskip 1cm X\in \uqn,\,Y\in \uqnm\leqno (1.4.7)$$
We can define a bilinear form $<\, ,\, >$ on
$\uq\times\uq$ by :
$$<X_1K_\lambda S(Y_1),\, Y_2K_\mu S(X_2)>=(X_1,\,Y_2)
(X_2,\,Y_1)q^{-(\lambda,\mu)/2}\leqno (1.4.8)$$\noindent
where $X_1,X_2\in \uqn,\,Y_1,Y_2\in \uqnm
,\,\lambda,\mu\in P$. This form is non degenerate.
\par There exists an algebra isomorphism from $\uqn$ to $\uqnm$ which maps
$E_i$
on $F_i$ for all $i$. Via this isomorphism, the restriction of the form 
$(\,,\,)$ on $\uqn\times\uqnm$ coincides with the one defined by Lusztig
in
[11, par 1].
\medskip\noindent
{\bf 1.5.} Define the maps :
$$\beta:\,\uqb\rightarrow\uqbm^*,\quad\beta(u)(v)=(u,\,v)\,
\,;$$  $$\zeta:\,\uq\rightarrow\uq^*,\quad\zeta(u)(v)
=<u,\,v>.$$\noindent
It follows from 1.4 that  : \bigskip\noindent
{\bf Lemma.} {\it With the previous notations, we have\par
\item{\rm (i)} $\beta,\,\zeta$ are injective.\par
\item{\rm (ii)} $\beta$ is an anti-homomorphism of
algebras.}\bigskip\noindent
Denote by $\rho$ the restriction homomorphism from $\uq^*$ onto
$\uqbm^*$. We know, see [2, Proposition 3.4], \medskip\noindent
{\bf Proposition.} {\it The restriction of  $\rho$ to $\cqgu$ is injective.
Moreover, for all $\lambda$ in $P^+$, we have
\smallskip\noindent
\item{(i)} For all  $e$ in $\uqn$,
$\rho(\zeta(eK_{-2\lambda}))=\beta(eK_{-\lambda})$.\par\noindent
\item{(ii)} there exists a (unique) subspace $E_\lambda$ of $\uqn$ such
that 
$\zeta(E_\lambda K_{-2\lambda})=C^+(\lambda)$.}\qed\bigskip\noindent
{\bf 1.6.} Let $u\mapsto \overline u$ be the $\K$-antihomomorphism of
$\uq$ such that $\overline E_i=E_i$, 
$\overline K_\lambda =K_{-\lambda}$,
$\overline F_i=F_i$. It is easily seen that 
$$<u,v>=(u,v)=(\overline u,\overline v)=<\overline u,\overline v>,
\,\,\; u\in\uqn,\, v\in\uqnm.$$
Let $\B$ be Lusztig's canonical basis of $\uqnm$, [11], which coincides
with 
Kashiwara's  global basis, [9].
Let $\B^*\subset\uqn$ be the dual basis in $\uqn$, i.e.
$(b^*,b')=\delta_{b,b'}$.
Let $\tilde E_i$, $\tilde F_i$ : $\uqnm\rightarrow\uqnm$ be the Kashiwara
operators, [{\it loc. cit.}]. For $b\in\B$, $\tilde E_i(b)$, resp.
$\tilde F_i(b)$, equals some $b'$ in $\B\cup\{0\}$ modulo
$q^{-1}\Z[q^{-1}]\B$.
The assignement $b\mapsto b'$ defines maps  $\tilde e_i$ and $\tilde f_i$ from
$\B$ to
$\B\cup\{0\}$.
For $b\in\B$, $1\leq i\leq n$, set $\varepsilon_i(b)=\hbox{Max}\{r,\,
\tilde e_i^r(b)\not=0\}$.\par
Let $L_i$, $1\leq i\leq n$, be the adjoint of the left multiplication
operator $F_i.\_$ for the form $(\,,\,)$ on $\uqn\times\uqnm$. Then,
$L_i$ is a quantum derivation of $\uqn$, [11, par 1] :
$$L_i(e_\alpha u)=L_i(e_\alpha)u+q^{(\alpha,\alpha_i)}e_\alpha L_i(u),
\;\; u\in\uqn,\, e_\alpha\in\uqn_\alpha, \alpha\in Q.$$
Set $L_i^{(r)}={1\over [r]_i!}L_i^r$.\medskip\noindent
 The following is a recollection of results about the canonical basis and its dual.
Assertions  (i) and (ii) can be read in 
 [11, 14.4.13,14.4.14]. Assertion (iii) is a standard consequence of 
 [11, 14.3.2 (c)] by dualization.
\medskip\noindent
{\bf Theorem 1.} {\it  For $b$, $b'$ in $\B$, we have :\par\noindent
\item {(i)} $bb'\in \Z[q,q^{-1}]\B$,  \par\noindent
\item {(ii)} $b^*b'^*\in \Z[q,q^{-1}]\B^*$, \par\noindent
\item {(iii)} $L_i^{(\varepsilon_i(b))}(b^*)=(\tilde
e_i^{\varepsilon_i(b)}(b))^*$,\par\noindent
\item {(iv)} $\overline b\in B$.}\qed\medskip\noindent
The following theorem states precisely the compatibility of the canonical
basis with
the finite dimensional highest weight modules. It can be found in
[9, Proposition 8.2].
\medskip\noindent
{\bf Theorem 2.} {\it  Fix $\lambda$ in $P^+$, 
$\lambda=\sum_i\lambda_i\varpi_i$. Set
$$\B_\lambda:=\{b\in \B,\;\varepsilon_i(\overline b)\leq \lambda_i,
1\leq i\leq n\}.$$\smallskip\noindent
Then $\B_\lambda=\{b\in \B,\; b.v_\lambda\not=0\}$ and
$\B_\lambda .v_\lambda$ is a basis of $\lql$.}\qed\medskip\noindent 
When $b$ belongs to $\B_\lambda$, and if no confusion occurs,
we shall use the same symbol for $b$ and $b.v_\lambda$, i.e. we set
$b=b.v_\lambda$.\medskip\noindent
{\bf 1.7.} Let $A=\C[q,q^{-1}]$. Let $\uanm$ be the $A$-submodule of 
$\uqnm$ generated by $\B$. Then, $\uanm$ is a free $A$-space and a
$A$-algebra.
Indeed, $\uanm$ is the $A$-algebra generated by the ${1\over
[m]_i!}F_i^m$,
see [8, Theorem 11.10 (b)]. Let $\ua$ be the sub-$A$-algebra
of $\uq$ generated by the ${1\over [m]_i!}F_i^m$, and the ${1\over
[m]_i!}E_i^m$.
Then, $\lal$ is a $\ua$-module.
Set $\lal=\uqnm.v_\lambda\subset\lql$. Then, $\lal$ is the $A$-space
generated
by $\B_\lambda$. By [8, Theorem 11.19], we know that $\B$ 
is compatible with specialisation :
$$\C\otimes_A\uanm\simeq U(\n^-),\;\; \C\otimes_A\lal\simeq V(\lambda),
\;\;\lambda\in P^+,$$ 
where $\C=A/(q-1)A$ as a $A$-module, $U(\n^-)$ is the (classical)
enveloping
$\C$-algebra of $\n^-$ and $V(\lambda)$ is the classical Weyl module with
highest weight
$\lambda$.  
\par\noindent
Let $\lal^*$ be the $A$-dual of $\lal$. Then, it has a natural
$\ua$-module structure
by $u.v^*(\_)=v^*(S(u).\_)$. The module $\lal^*$ specializes at $q=1$
onto
the dual $\g$-module $V(\lambda)^*$.\medskip\noindent
{\bf 1.8.} Fix $\lambda$ in $P^+$ and $w$ in $W$. 
We know that $\lql$ verifies the Weyl character formula ;
we denote by $v_{w\lambda}$ an extremal vector  
of weight $w\lambda$. Then, the $\uqb$-module $V_{q,w}(\lambda):=
\uqn v_{w\lambda}$ verifies
the Demazure character formula. We know, [9, Theorem 12.4], [12,
5.3-5.4], that 
\medskip\noindent
{\bf Theorem.} {\it There exists
a subset $\B_w$ of $\B$ such that $V_{q,w}(\lambda)$ is spanned
by
$\B_w.v_\lambda$. Moreover, if $b$ is in $\B_w$ then $\Delta(b)
\in <\B_w>\otimes <K_\mu .\B_w, \mu\in P>$.}\qed\medskip\noindent
In particular, the orthogonal $V_{q,w}(\lambda)^\bot$ in
$\lql$ of the Demazure module $V_{q,w}(\lambda)$ is generated as a space
by
$\B(\lambda)\cap \B\backslash\B_w$ and the dual $V_{q,w}(\lambda)^*$
is generated by the image of $\B(\lambda)\cap\B_w$ by the quotient
morphism.
\par
As in 1.7, this allows us to define $A$-forms for Demazure modules. We
denote by
$V_{A,w}(\lambda)$ the $A$-space generated by $\B(\lambda)\cap \B_w$.
It specializes for $q=1$ to the classical Demazure module $V_w(\lambda)$.
\bigskip\noindent
{\bf 2. A multiplicative property and Littelmann's parametrization of the
dual canonical basis.}\bigskip\noindent
{\bf 2.1.} For $\lambda$ in $P^+$ and $b$ in $\B_\lambda$, let
$b_\lambda^*$
be the element of $\lql^*$ such that
$b_\lambda^*(b'.v_\lambda)=\delta_{b,b'}$,
$b'\in \B_\lambda$, where $\delta$ means the Kroenecker symbol. 
\medskip\noindent
{\bf Lemma.} {\it For all $\lambda$ in $P^+$ and $b$ in $\B_\lambda$,
then $\zeta(b^*K_{-2\lambda})=b_\lambda^*\otimes v_\lambda$,
$\beta(b^*K_{-\lambda})=\rho(b_\lambda^*\otimes
v_\lambda)$.\smallskip\noindent
Proof.} By Lemma 1.5, we only need to prove that 
$<b^*K_{-2\lambda},-->=b_\lambda^*\otimes v_\lambda$. As $\uq=\uqbm\oplus
(\uqbm\otimes\uqn\cap\Ker(\varepsilon))$, we only need to prove this on
$\uqbm$, whose basis is given by $(b'K_\mu, \;b'\in\B,\mu\in P)$.
By (1.4.7) and (1.4.8)
$$<b^*K_{-2\lambda},b'K_\mu>=(b^*,b')q^{\lambda,\mu}
=\delta_{b,b'}q^{\lambda,\mu}=b_\lambda^*(b'K_\mu. v_\lambda).$$
This implies the lemma.\qed\medskip\noindent
{\bf Remark.} By Proposition 1.5, the lemma implies that 
$E_\lambda$ is spanned  by $\B_\lambda^*$.\medskip\noindent
Let $\aqgu$ be the sub-$A$-module of $\cqg$ generated by the 
$b_\lambda^*\otimes v_\lambda$, $\lambda\in P^+$ and $b\in \B_\lambda$.
Let $d_{b,b'}^{b''}$ be the coefficient of $b''^*$ in the product 
$b^*b'^*$.\medskip\noindent
{\bf Proposition.} {\it We have \par\noindent
\item{(i)} $\K\otimes_A \aqgu=\cqg$\par\noindent 
\item{(ii)} $\aqgu$ is a $A$-algebra.  \par\noindent
\item{(iii)} if $d_{b',b}^{b''}$ in non zero, then $b\in \B_\lambda$,
$\b'\in\B_{\lambda'}$ implies $b''\in
\B_{\lambda+ \lambda'}$, \par\noindent
\item{(iv)} $\aqgu/(q-1)\aqgu\simeq\C[G/U]$, where $\cgu$ 
is the 
$\C$-algebra of regular functions on $G/U$.\smallskip\noindent
Proof.} The $A$-basis $b_\lambda^*$ of $\lal^*$ is a 
$\K$-basis of $\lql^*$. Hence, (i) is clear.\par\noindent
From Lemma 1.5, Proposition 1.5, and the previous lemma, we have
$$\beta^{-1}[(b_\lambda^*\otimes v_\lambda).({b'_{\lambda'}}^*\otimes
v_{\lambda'})]=\beta^{-1}({b'_{\lambda'}}^*\otimes
v_{\lambda'}).\beta^{-1}(b_{\lambda}^*\otimes v_{\lambda})$$
$$=(b'^*K_{-\lambda'}).(b^*K_{-\lambda})=q^{-(\lambda',\omega(b^*))}
b'^*b^*K_{-\lambda-\lambda'}=q^{-(\lambda',\omega(b^*))}(\sum
d_{b',b}^{b''}{b''}^*)K_{-\lambda-\lambda'}.$$
We know that $(b_\lambda^*\otimes v_\lambda).({b'_{\lambda'}}^*\otimes
v_{\lambda'})$ belongs to $C^+(\lambda+\lambda')$ which is generated by 
$b^0\otimes  v_{\lambda+\lambda'}$,  $b^0\in \B_{\lambda+\lambda'}$. 
Hence, applying $\beta^{-1}$, we obtain (iii) and  
$$(b_\lambda^*\otimes v_\lambda).({b'_{\lambda'}}^*\otimes v_{\lambda'})
=q^{-(\lambda',\omega(b^*))}\sum
d_{b',b}^{b''}({b_{\lambda+\lambda'}''}^*\otimes
v_{\lambda+\lambda'}).\leqno
(2.1.1)$$
This gives (ii). Now, (iv) is clear by 1.7 and the fact that
specialization commutes with tensor product.\qed\medskip\noindent
{\bf Corollary.} {\it Fix $b$, $b'$, $b''$ in $\B$ with
$d_{b,b'}^{b''}$ non zero. Then, for all $i$, $1\leq i\leq n$, we have
$\varepsilon_i(b'')\leq
\varepsilon_i(b)+\varepsilon_i(b')$.\smallskip\noindent
Proof.} By applying the antiautomorphism ${ }^{--}$, we obtain 
 that $d_{\overline b',\overline b}^{\overline b''}
=d_{b,b'}^{b''}$. Set $\lambda=\sum_i\varepsilon_i(b)\varpi_i$ and
$\lambda'=\sum_i\varepsilon_i(b')\varpi_i$. As ${ }^{--}$ is an
involution, we
deduce from  [1.6, Theorem 2] that $\overline b\in \B_\lambda$ and
$\overline b'\in \B_{\lambda'}$. This gives 
$\overline b''\in \B_{\lambda+\lambda'}$ by (iii). Hence, we obtain the
corollary 
from [1.6, Theorem 2].\qed\medskip\noindent
{\bf Remark.} In the simply laced case,
this corollary is easily obtained by the positivity property of the dual 
canonical basis, i.e. $d_{b,b'}^{b''}\in \N[q,q^{-1}]$ by [11].
For general $\g$, we can conclude by a standard argument, transmitted by
E. Vasserot. It is based on the realization of 
$d_{b,b'}^{b''}$ in terms of traces of an automorphism of a diagram on
spaces arising
from perverse sheaves, see [11, 14.4.14].\medskip\noindent
{\bf 2.2.} We introduce Littelmann's parametrization of the 
(dual) canonical basis. Fix a reduced decomposition of the longest element
of the Weyl group
$w_0$ : $\tilde w_0=s_{i_1}\ldots s_{i_N}$, where $N=\dim\n$.
For all $u$ in $\uqn$ and $1\leq i\leq n$, set
$$a_i(u)=\hbox{Max}\{r,\,
L_i^r(u)\not=0\},\;\;\Lambda_i(u)=L_i^{(a_i(u))}(u).$$
For all $b$ in $\B$, set
$$A_{\tilde w_0}(b)=(a_{i_N}(b^*),a_{i_{N-1}}(\Lambda_{i_N}(b^*)),
\ldots,a_{i_1}(\Lambda_{i_2}\ldots\Lambda_{i_N}(b^*)))\in \N^N.$$
This parametrization can be found in [10, par 1]. It coincides with the
parametrization
in [1, 3.2] by [1.6, Theorem 1 (iii)].
We now present a theorem due to Littelmann, [10, par 1], see also
[1, 3.10].
\medskip\noindent
{\bf Theorem.} {\it The map  $A_{\tilde w_0}$ embeds $\B$ into $\N^N$.
Let $\co$ be its image. Then, $\co$ is the set of integral points of a rational 
convex polyhedral cone of $\R^N$.
Moreover, set $\gwo:=\{(\lambda,\psi)\in P^+\times\co,\;\psi\in
A_{\tilde w_0}(\B_\lambda)\}$. Then, $\gwo$ is the set of integral points
of a rational
convex polyhedral cone of $\R^{n+N}$.}\qed\medskip\noindent
Note that equations of this cone can be given in [1, 3.10].
For all $\psi$ in $\co$, set 
$b^\psi=b_{\tilde w_0}^\psi=A_{\tilde w_0}^{-1}(\psi)$.
We also set $d_{\psi,\psi'}^{\psi''}:=d_{b^\psi,b^{\psi'}}^{b^{\psi''}}$.
\medskip\noindent
{\bf 2.3.} Let $\prec$ be the lexicographical ordering of $\N^N$.
We have\medskip\noindent
{\bf Theorem.} {\it Fix a reduced decomposition $\tilde w_0$ of the
longest 
element of the Weyl group. Let $b$, $b'$, $b''$ be in $\B$, 
$A_{\tilde w_0}(b)=\psi$, $A_{\tilde w_0}(b')=\psi'$
Then, $d_{b,b'}^{b''}$ non zero implies $A_{\tilde w_0}(b'')
\prec \psi+\psi'$. Moreover, 
$d_{\psi,\psi'}^{\psi+\psi'}$ is a power of $q$.\smallskip\noindent
Proof.} First remark that $a_i(uv)=a_i(u)+a_i(v)$, $u$, $v$ in $\uqn$, 
by the quantized Leibniz rule. Write
$$b^*b'^*=\sum d_{b,b'}^{b''}{b''}^*\leqno (*)$$ 
Let $\psi''=A_{\tilde w_0}(b'')$, for $b''$ in
the sum. Then, $\psi''_1\leq \psi_1+\psi'_1$ by Corollary 2.1
and [1.6, Theorem 1 (iii)]. This gives the first step of the induction.
Now, applying $L_{i_N}^{\psi_1+\psi'_1}$ and the Leibniz rule gives 
$\xi_1\Lambda_{i_N}(b^*)\Lambda_{i_N}(b'^*)=\sum d_{b,b'}^{b''}
\Lambda_{i_N}(b''^*)$, where the sum is taken over the elements
$b''$ such that $\psi''_1= \psi_1+\psi'_1$, and where $\xi_1$ is a 
power of $q$. Now remark that, by 
[1.6, Theorem 1 (iii)],
$\Lambda_{i_N}(b^*)$, $\Lambda_{i_N}(b'^*)$, and the
$\Lambda_{i_N}(b''^*)$ all belong to the dual of the canonical basis. 
So, we can proceed as for the first step by induction. We then have the
first
assertion of the theorem. For conclusion, note that we obtain
$$\xi_N\Lambda_{i_1}\ldots\Lambda_{i_N}(b^*)
\Lambda_{i_1}\ldots\Lambda_{i_N}(b'^*)=\sum d_{b,b'}^{b''}
\Lambda_{i_1}\ldots\Lambda_{i_N}(b''^*),$$
where the sum is taken over the elements
$b''$ such that $\psi''_k= \psi_k+\psi'_k$, for all $k$
 and where $\xi_N$ is a 
power of $q$.  Hence, by Theorem 2.2 there is at most one element in the
sum
and for this element we have $A_{\tilde w_0}(b'')=\psi+\psi'$. By 
[10, par 1],\medskip\noindent
{\bf Lemma.} {\it $\Lambda_{i_1}\ldots\Lambda_{i_N}(B^*)=1$ for all
elements
$B$ of the canonical basis.\qed\medskip\noindent}
Hence the coefficient $d_{\psi,\psi'}^{\psi+\psi'}$ is a power of $q$.
This finishes the proof of the theorem.
\qed\medskip\noindent
{\bf 2.4.} By 1.8, the results of 2.1, 2.2 can be easily generalized 
to quotients  of $\cqgu$ which correspond to Demazure modules. Indeed, fix
$w$ in $W$ and
let $\B_{\overline w}$ be the complement of $\B_w$ in $\B$. Set
$$I_{A,w}:=\bigoplus_{\lambda\in P^+,\,
b\in \B_{\overline w}\cap\B(\lambda)}Ab^*\otimes v_\lambda
=\bigoplus_\lambda V_{A,w}(\lambda)^\bot\otimes v_\lambda.$$
Then, $I_{A,w}$ is the orthogonal $<\B_w>^\bot$ of $<\B_w>$ in $\aqgu$. 
By Theorem 1.8, $I_{A,w}$ is an ideal of $\aqgu$. 
We have the following decomposition 
for the quotient algebra : $$\aqgu/I_{A,w}=
\bigoplus_\lambda V_{A,w}(\lambda)^*\otimes v_\lambda.$$
From [10, par 1], we have :
\medskip\noindent
{\bf Theorem.} {\it Let $\tilde w=s_{i_1}\ldots s_{i_p}$ 
be a reduced decomposition of $w$. Then there exists a reduced
decomposition
 $\tilde w_0=s_{i_N}\ldots s_{i_1}$ of $w_0$. For this decomposition, we
 have $A_{\tilde
w_0}(B_w)=\co\cap(\N^p\times\{0\}^{N-p})$.}\qed\medskip\noindent
{\bf Definition.} For all $w$, a reduced decomposition $\tilde
w_0=s_{i_N}\ldots s_{i_1}$ of $w_0$ such that $w=s_{i_1}\ldots s_{i_p}$
will be called adapted to $w$.\bigskip\noindent 
{\bf 3. Specialization.}\bigskip\noindent
\def\pil{\phi^{k,l}}\def\psii{\psi^{k}}
{\bf 3.1.} Fix a reduced decomposition
$\tilde w_0$ of $w_0$. At this stage of the article, we can
construct 
a
$\N^N$-filtration of the algebra $\cgu$ such that the associated graded
algebra is the algebra of the semigroup $\gwo$ 
\def\blp{b_{\lambda,\psi}^*}\def\blph{b_{\lambda,\phi}^*} 
To be more precise, let $\blp$ in $\cgu$ be the image of the 
element $(b_\lambda^\psi)^*\otimes v_\lambda$ by the  morphism of
specialization at
$q=1$, see Proposition 2.1 (iv). We have by Theorem 2.3 and (2.1.1)
:\medskip\noindent 
{\bf Proposition.} {\it The spaces $\F_\psi:=< \blph,\;
 (\lambda,\phi)\in\gwo,\;\phi\prec\psi\}$, $\psi\in\co$ define a
filtration 
 of the algebra $\cgu$. The graded associated algebra is naturally
isomorphic to
 the $\C$-algebra of the semigroup $\gwo$.}\qed\medskip\noindent
{\bf 3.2.} What results from Proposition 3.1 is that there exists a
finite sequence of degenerations of the flag variety which ends into a
toric variety but what we want is a "degeneration in one step". Hence, we
need
a linear form  $\N^N\rightarrow \N$, which satisfies some strict
inequalities, and which transforms the $\N^N$-filtration of $\cgu$ into a
$\N$-filtration. 
This is made possible because the cone $\gwo$ is convex polyhedral and
hence has a finite presentation. We start by a lemma.\medskip\noindent
{\bf Lemma.} {\it Let $\{\psi^1,\ldots,\psi^K\}$ be a finite set of points
in $\N^N$, and for all $k$, $1\leq k\leq K$, let
$\phi^{k,l}$, $1\leq l\leq K_k$, be a finite number of points in $\N^N$
such that $\pil\prec\psii$ for all $l$. 
Then, there exists a linear form $e$ : $\N^N\rightarrow\N$ such that
$e(\pil)<e(\psii)$ for all $k$ and all $l$.\smallskip\noindent
Proof.} Let $a_k$, $1\leq k\leq N$ be the linear form of $\Q^N$ which maps
an element of $\Q^N$ to its $k$-th coordinate. Set $I_0=\{(\pil,\psii),\;
1\leq k\leq K,\,1\leq l\leq K_k\}$.
For $s$, $1\leq s\leq N-1$, set $I_s=\{(\pil,\psii)\in I_0,\; (\pil)_m=
(\psii)_m, 1\leq m\leq s\}$. All these sets are finite. \par
As a first step of our induction, define the linear form $e_N=a_N$ and fix
$\varepsilon_N$
in $\Q^+$ such that $e_N(\pil)\varepsilon_N<(\psii)_{N-1}-(\pil)_{N-1}$,
for all $(\pil,\psii)$ in $I_{N-2}\backslash I_{N-1}$. Define a linear
form by $e_{N-1}=a_{N-1}+\varepsilon_Ne_N$. By construction, we have 
$e_{N-1}(\pil)<e_{N-1}(\psii)$ for all $(\pil,\psii)$ in $I_{N-2}$.
Fix now $\varepsilon_{N-1}$
in $\Q^+$ such that $e_{N-1}(\pil)\varepsilon_N<
(\psii)_{N-2}-(\pil)_{N-2}$,
for all $(\pil,\psii)$ in $I_{N-3}\backslash I_{N-2}$. 
Define a linear form by $e_{N-2}=a_{N-2}+\varepsilon_{N-1}e_{N-1}$. By
construction and induction we have $e_{N-2}(\pil)<e_{N-2}(\psii)$ for all
$(\pil,\psii)$ in $I_{N-3}$. By induction, we obtain a form 
$e_1$ which verifies $e_{1}(\pil)<e_{1}(\psii)$ for all $(\pil,\psii)$ in
$I_{0}$. We can multiply $e_1$ by a positive integer in order to obtain a
$\N$-form which verifies the lemma.\qed\medskip\noindent 
Let $(\lambda_i,\psi_i)$, $1\leq i\leq K$, be the minimal set of
generators of
$\gwo$ and
$(\sum_in_i^k(\lambda_i,\psi_i)-\sum_im_i^k(\lambda_i,\psi_i))$,
$n_i^k$, $m_i^k\in\N,$ be the (finite) minimal set of generators of the 
relations. Then, from Proposition 3.1, we obtain :\medskip\noindent
{\bf Proposition.} {\it The (commutative) algebra $\cgu$ is defined by generators
 $b_{\lambda_i,\psi_i}^*$, $1\leq i\leq K$, and relations 
 $$\prod b_{\lambda_i,\psi_i}^{*n_i^k}=\prod b_{\lambda_i,\psi_i}^{*m_i^k}+
 \sum\C b_{\mu,\phi}\leqno (3.2.1)$$ for 
 $\phi\prec\sum_in_i^k(\lambda_i,\psi_i)=\sum_im_i^k(\lambda_i,\psi_i)$.
 \smallskip\noindent
Proof.} The algebra $\cgu$ has a natural $P^+$-grading defined by
$\cgu_\lambda=V(\lambda)^*\otimes v_\lambda$. By construction, 
$\F_\psi=\oplus_\lambda\F_{\psi,\lambda}$, where $\F_{\psi,\lambda}=
\F_\psi\cap\cgu_\lambda$. Clearly, $\F_{\psi,\lambda}$ is finite
dimensional. Hence, by induction, we obtain that generators of the 
$\co$-graded algebra of Proposition 3.1 lift to generators of the
algebra $\cgu$. Hence, we have the first part of the proposition. 
\par Now, Proposition 3.1 implies the relations (3.2.1). Let
$\C[X_i, 1\leq i\leq K]$ be the polynomial algebra with $K$
indeterminates.
There exists a surjective morphism $\varphi$ : $\C[X_i, 1\leq i\leq
K]/J\rightarrow
\cgu$, where $J$ is the ideal generated by the relations resulting from
(3.2.1). This morphism maps $X_i$ to  $b_{\lambda_i,\psi_i}^*$ for all $i$. Endow
$\C[X_i, 1\leq i\leq
K]/J$ with the quotient filtration $<\prod X_i^{m_i}+J,\; \sum m_i\psi_i\prec 
\psi>_{\psi\in\co}$. Then, the associated graded algebra is defined by generators
Gr$X_i$, $1\leq i\leq K$, and relations 
$\prod \Gr \,X_i^{n_i^k}=\prod \Gr \,X_i^{m_i^k}$. Now, endow $\cgu$ with its 
$\co$-filtration, see 3.1. By Proposition 3.1,
Gr$\varphi$ is an isomorphism from 
Gr$\C[X_i, 1\leq i\leq K]/J$ onto Gr$\cgu$.
This implies that
$\varphi$ is an isomorphism. This finishes the proof of the
proposition.\qed\medskip\noindent
Let $\psi^k:=\sum_in_i^k\psi_i$, and for each $k$, let the $\phi^{k,l}$ 
 be the 
 elements $\phi$ occuring in (3.2.1). 
 Fix a $\N$-form $e$ as in Lemma 3.2.  Then,
 \medskip\noindent
{\bf Corollary.} {\it The spaces $\F_m(\cgu)=< \prod
b_{\lambda_i,\psi_i}^{s_i^k},
\;e(\sum_is_i^k\psi_i)\leq m>$, $m\in\N$, define a filtration of $\cgu$.
The 
graded associated algebra is $\C[\gwo]$.}\qed\medskip\noindent 
{\bf 3.3.} We now give some analogous results for Demazure modules. Fix $w$
in $W$ with length $p$. 
Set $\gw:=A_{\tilde w_0}(\B_w)$. Then,\medskip\noindent
{\bf Lemma.} {\it For all reduced decompositions $\tilde w_0$ of $w_0$, the
set $\gw$ is a finite union of rational convex polyhedral cones. Moreover, if 
$\tilde w_0$ is adapted to $w$, see Definition 2.4, then $\gw$ is a
$p$-dimensional face of $\gwo$. In particular, it is a rational convex polyhedral
cone.
\smallskip\noindent
Proof.} Fix $w$ in $W$ and fix a reduced decomposition $\tilde w_0$
adapted to $w$. By Theorem 2.4,
$\gw$ is a $p$-dimensional face of $\gwo$ and so it is a convex polyhedral
cone. Let $\tilde w_0'$ be another reduced decomposition of $w_0$ which is
not supposed adapted to $w$. Then, by [1, 3.3], $\gwp$ is the image
of
$\gw$ by a continuous piecewise linear map $l$. To be more precise, there exists a finite set
of (convex) cones $\co_i$ in $\R^{n+N}$ such that 
$\bigcup \co_i=\R_{\geq 0}^{n+N}$ and such that $l$ is linear on $\gw\cap\co_i$.
 Hence, we have
the lemma.\qed\medskip\noindent 
 We still denote by $v_\lambda$ a highest weight vector of the classical Weyl
module
 $V(\lambda)$. The algebra $\aqgu$ specializes for $q=1$ 
 onto the algebra 
$\cgu=\bigoplus_\lambda V(\lambda)^*\otimes v_\lambda.$
The algebra $\aqgu/I_{A,w}$ specializes for $q=1$ onto the quotient
algebra 
$\cgu/I_w=\bigoplus_\lambda V_w(\lambda)^*\otimes v_\lambda$, where
$I_w:=\bigoplus_{\lambda\in P^+} 
V_w(\lambda)^\bot\otimes v_\lambda.$ \par 
Let  
$\tilde w_0$ be adapted to $w$. By the previous lemma, $\gw$ is
generated as a semigroup by a part of the minimal set of generators
$(\lambda_i,\psi_i)$ of $\gwo$ and with the corresponding relations.
This implies :\medskip\noindent
{\bf Theorem.} {\it Choose a reduced decomposition $\tilde w_0$  adapted
to $w$. Then, the graded associated algebra of the quotient
filtration of $(\F_m)_{m\in\N}$ on $\cgu/I_w$ is the $\C$-algebra 
$\C[\gw]$ of the semigroup $\gw$.}\qed\medskip\noindent
{\bf Remark.} Note that this filtration (and hence the corresponding
degeneration) depends on the choice of the
$\N$-form $e$ of Lemma 3.2, and of the reduced decomposition $\tilde w_0$
adapted to $w$. For a general reduced decomposition $\tilde w_0'$ of
$w_0$, Lemma 3.3 implies that the spectrum of the associated graded
algebra is a union of irreducible components which are toric varieties.    

\vskip 1 cm \noindent{\petittitre\centerline{ACKNOWLEDGMENTS}}\vskip 1cm
We wish to thank R. Yu for a general survey of the Degeneration Problem. This work owes
much to helpful conversations with M. Brion and O. Mathieu. We are
grateful to P. Baumann for true criticisms and to J. Germoni for
encouragements and fan constructions. 
\vskip 1 cm \noindent{\petittitre\centerline{BIBLIOGRAPHY}}\vskip 1cm
\parindent=0cm
[1]. A. BERENSTEIN and A. ZELEVINSKY. {\it Tensor product
multiplicities, Canonical bases
and Totally positive varieties}, ArXiv:Math.RT/9912012.\smallskip 
[2]. P. CALDERO. {\it G\'en\'erateurs du centre de $U_q(sl(N+1))$},
Bull.
Sci. Math., 118, (1994),
177-208.\smallskip
[3]. P. CALDERO. {\it Elements ad-finis de certains groupes
quantiques},
 C.R.Acad.
Sci. Paris, t. 316, Serie I, (1993), 327-329.\smallskip
[4]. R. CHIRIVI. {\it LS Algebras and applications to Schubert Varieties},
Trans. Groups, Vol. 5, No 3, 245-264, (2000).\smallskip
[5]. R. DEHY.
{\it Polytopes associated to Demazure modules of 
symmetrizable Kac-Moody algebras of rank two}, 
J. Algebra, 228, No.1, 60-90, (2000).\smallskip
[6]. R. DEHY, R.W.T. YU.
{\it Lattice polytopes associated to certain Demazure modules of
$sl_{n+1}$},
J. Algebr. Comb., 10, No.2, 149-172 (1999).\smallskip
[7]. N. GONCIULEA, V. LAKSHMIBAI.
{\it Degenerations of flag and Schubert varieties to toric varieties}, 
Transform. Groups 1, No.3, 215-248 (1996).\smallskip 
[8]. J. C. JANTZEN,  
{\it Lectures on quantum groups}, 
Graduate Studies in Mathematics, 6, Providence, 
 American Mathematical Society. vii, 266 p, (1996).\smallskip
[9]. M. KASHIWARA. {\it On Crystal Bases}, Canad. Math. Soc.,
Conference Proceed., 16, (1995), 155-195.\smallskip
[10]. P. LITTELMANN, 
{\it Cones, crystals, and patterns}, 
Transform. Groups 3, No.2, 145-179 (1998).\smallskip
[11].  G. LUSZTIG. {\it Introduction to quantum groups}, Progress in
Mathematics,
 110, Birkauser, (1993).
\smallskip 
[12]. G. LUSZTIG.
{\it Problems on canonical bases},
Haboush, William J. (ed.) et al., Algebraic groups and their 
generalizations: quantum and infinite-dimensional methods,
 Proc. Symp. Pure Math. 56, Pt. 2, 169-176 (1994).\smallskip
[13]. M. REINEKE. {\it Multiplicative Properties of Dual Canonical Bases
of
Quantum Groups}, J. Alg., 211, (1999), 134-149.\smallskip
[14]. M. ROSSO. {\it Analogues de la forme de Killing et du th\'eor\`eme
 de Harish-Chandra
pour les groupes quantiques}, Ann. Sci. Ec. Norm. Sup., 23, (1990),
445-467.\smallskip
[15]. T. TANISAKI. {\it Killing forms, Harish-Chandra isomorphisms, and
universal
{\cal R}-matrices for quantum algebras}, Int. J. Mod. Phys. A, Vol. 7,
 Suppl. 1B,
 (1992),  941-961.\smallskip

{\it Institut Girard Desargues, UPRES-A-5028  \par
Universit\'e Claude Bernard Lyon I, Bat 101\par
69622 Villeurbanne Cedex, France}\smallskip\noindent
e-mail : caldero@desargues.univ-lyon1.fr
\medskip\noindent 
Keywords : Toric Degeneration, Schubert Varieties, Canonical Basis.
\end